\documentclass[11pt]{article}

\usepackage{graphicx,amsmath,amssymb,amsthm,subfigure,url,xspace}
\usepackage{fullpage,amsbsy,equation,enumerat}
\usepackage{bbm,color}

\usepackage{pgfplots}

\usetikzlibrary{arrows.meta}
\usetikzlibrary{shapes,backgrounds}
\usetikzlibrary{patterns}
\usetikzlibrary{decorations.markings}

\pgfplotsset{compat = newest}

\pgfplotsset{
  every axis/.append style={
      axis x line=middle,    
       axis y line=middle,    
       axis line style={<->,color=blue}, 
       xlabel={$x$},          
       ylabel={$y$},          
       }
 }

\newcommand{\mb}[1]{{\boldsymbol #1}}
\newcommand{\mq}{\mb{q}}
\newcommand{\dt}{\, \hbox{d}t}
\newcommand{\dx}{\, \hbox{d}\mb{x}}
\newcommand{\R}{{\mathbb R}}
\newcommand{\1}{{\mathbbm{1}}}
\renewcommand{\div}[1]{\hbox{div} \, #1}
\newcommand{\blot}[1]{}

\newtheorem{theo}{Theorem}
\newtheorem{lemma}[theo]{Lemma}
\newtheorem{cor}[theo]{Corollary}
\newtheorem{prop}[theo]{Proposition}

\makeatletter
\def\ltag#1{%
\stepcounter{equation}%
\tag*{}%
\def\df@tag{(\theequation)\llap{\rlap{#1}\hspace{\columnwidth}}}%
}
\makeatother

\numberwithin{equation}{section}

\title{The Least Action Admissibility Principle}
\author{H.~Gimperlein\thanks{Engineering Mathematics, University of Innsbruck, Innsbruck, Austria} \and M.~Grinfeld\thanks{Department of Mathematics and Statistics, University of Strathclyde, Glasgow, G1 1XH, UK} \and R.~J.~Knops\thanks{The Maxwell Institute of Mathematical Sciences and School of Mathematical and Computing Sciences, Heriot-Watt University, Edinburgh, EH14 4AS, Scotland, UK} \and M.~Slemrod\thanks{Department of Mathematics, University of Wisconsin, Madison, WI 53706, USA}}

\begin{document}
\maketitle

\begin{abstract}
\noindent  This paper provides a new admissibility criterion for choosing physically relevant weak solutions of the equations of Lagrangian and continuum mechanics when non-uniqueness of solutions to the initial value problem occurs. The criterion is motivated by the classical least action principle but is now applied to initial value problems which exhibit non-unique solutions. Examples are provided to Lagrangian mechanics and the Euler equations of barotropic fluid mechanics. In particular, we show the least action admissibility principle prefers the classical two shock solution to the Riemann initial value problem to certain solutions generated by convex integration. {On the other hand, Dafermos's entropy criterion prefers convex integration solutions to the two shock solutions. Furthermore, when the pressure is given by $p(\rho)=\rho^2$, we show that the two shock solution is always preferred whenever the convex integration solutions are defined for the same initial data.}
\end{abstract}
\section*{Introduction}

In his beautiful monograph \cite[pp. 119--120]{Ekeland}, Ivar Ekeland
has written: ``Think of a massive body, a small one like an electron
or a large one like a billiard ball, starting from $A$ and ending up
at $B$. What path will it take? The answer we get from classical
physics, in the absence of any external force, is a straight line.

Feynman's answer \cite{Feynman} is every path from $A$ to $B$ is
possible, from the straight one to the most crooked one you can
imagine. To find how likely a given path is, one has to compute the
action (yes, the classical action as defined by Maupertuis, Euler,
Lagrange, Hamilton, Jacobi, the old crowd) along that path \ldots\, If
one gets into the mathematics, one finds that the [most likely paths]
are the ones that make  the action stationary.''


In the language of the mathematical theory of conservation laws, we
might say Feynman introduces the idea of a selection or admissibility
criterion to the non-uniqueness of particle trajectories in quantum
mechanics, in fact the classical least action principle. The goal of
this work is to follow in Feynman's footsteps and
introduce the least action principle as an admissibility criterion
when confronted by the situation when the Euler equations of classical
compressible fluid mechanics  possess non-unique solutions to
initial value problems.

Here we define the Lagrangian $L$ of a motion (particle or fluid) as
the pointwise difference of the kinetic and potential energy. (Brenier
\cite{Brenier} {assigns} the potential energy to have the
opposite sign and hence his definition of the Lagrangian has the sum
of the kinetic and potential energies.) The action for particle motion
is given by
\[
  A(\mq) = \int_{t_0}^{t_1} L(\mq,\dot{\mq}) \dt, \; \; t_1 > t_0,
\]
where $\mq, \dot{\mq}$ denote the position and the velocity of the
particle at time $t$ which was initially at $\mq(t_0) = \mq_0$,
$\dot{\mq}(t_0) = \mq_1$.

The Lagrangian for the motion of a compressible barotropic fluid is
given by
\[
  L(\rho, \mb{v}) = \frac12 \rho |\mb{v}|^2 - \rho \epsilon(\rho).
\]
Here $\rho$ is the mass density, $\mb{v}$ is the velocity, and $\epsilon(\rho)$ is the specific internal energy,
$\epsilon'(\rho)=p(\rho)/\rho^2$, with $p$ the pressure and $p'(\rho) > 0$ throughout this paper. The action of
a fluid motion confined to a domain $\Omega$  is given by
\[
  A(\rho, \mb{v}) = \int_{t_0}^{t_1} \int_\Omega
  L(\rho, \mb{v}) \, \hbox{d}\mb{x} \dt,  \; \; t_1 > t_0.
\]
Notice we use the Eulerian description for fluid motion as this is the
most convenient for our examples.

 The \emph{least action admissibility principle} (LAAP) {states} that an ``entropic'' weak solution to the relevant particle
or fluid motion is admissible if the action for this motion is less than or equal to the action of all
other ``entropic'' motions with same initial and boundary data. Here ``entropic'' means the weak solutions considered satisfy the usual energy-entropy conditions as given in Section 3.
Thus LAAP is viewed as a supplement and not as a replacement to the usual entropy criterion.

Notice the difference between the classical least action principle and LAAP. In the classical least action principle we are given the Lagrangian, initial and terminal positions, i.e.~two pieces of information beyond the Lagrangian. With this formulation, stationarity of the action will produce the Euler-Lagrange equations. In LAAP {the two pieces of data are} initial position and momentum. There is no need to look for stationarity of the action, since the Lagrangian is derived from the already known balance laws, {the action now serving the new purpose as an admissibility criterion.} In practice, while it would be desirable to compare one weak solution to an initial value problem with all others, here we consider the more accessible task of comparing one known weak solution with other known weak solutions. {Consequently,} in this article we know explicit formulas for the action.


The use of the least action principle to prove existence and
uniqueness for the incompressible Euler equations originates in the
paper of Ebin and Marsden \cite{EbinMarsden} and was developed by Brenier \cite{Brenier}. These results were
improved in the papers of Shnirelman \cite{Shnirel1, Shnirel2} and the
subject has been surveyed by Brenier in \cite{Brenier2}. In their paper on Rayleigh-Taylor mixing Gebhart et al.~\cite{Gebhart} have used a least action principle as a selection mechanism for subsolutions necessary for the convex integration process. We also note the more recent work on dual variational least action approaches to existence of solutions to the incompressible fluid equations as given in \cite{aa24} and the references included there. To our
knowledge, the use of the least action principle as an admissibility
criterion for weak solutions of the compressible Euler equations is
new.

The major motivation for this work begins with the results of De
Lellis and Sz\'ekelyhidi \cite{deLelSz1, deLelSz2}, {who use} the
method of convex integration to construct infinitely many solutions to
both incompressible and compressible Euler systems, all emanating from
the same initial data and satisfying the standard energy--entropy
admissibility criterion. A complete list of contributors to this
project is found in the recent monograph of S.~Markfelder
\cite{Markfelder}.

As noted above, since the weak solutions constructed by the convex
integration method satisfy the usual energy--entropy admissibility
criterion {but are non-unique}, the next step was to use a more precise criterion, i.e. the
entropy rate criterion of C. ~M.~Dafermos \cite{Dafermos}. In his paper
\cite{Feireisl}{,} E. Feireisl showed that under some mild restrictions,
the entropy rate criterion rules out all weak solutions to the
barotropic compressible Euler equations constructed by convex
integration. The issue might be considered settled if not for
the surprising result of Chiodaroli and Kreml \cite{ChioKr}, {who found a two-dimensional Riemann problem in which}
the convex integration solutions could be preferred, according to the
entropy rate criterion, to the perhaps physically expected two shock
solution. Hence it was strong motivation to consider the role of the
least action admissibility principle in this example. Here we give
sufficient conditions on the Riemann data for which the two shock
solution will be preferred according to the least action
admissibility criterion.

This paper is divided into five sections after this introduction. Section 1 provides a physical interpretation of LAAP. More precisely, it motivates LAAP as a criterion for choosing the material motion that minimizes cost or effort. Section 2 considers an oscillator {treated} by Dafermos \cite{Dafermos2}. The example is an elementary system of ordinary differential equations exhibiting non-uniqueness for the initial value problem. Dafermos has shown that uniqueness may be recovered via application of vanishing viscosity and entropy rate criteria. {The} same uniqueness result follows from LAAP. Section 3 {treats} non-uniqueness results of Akramov and Wiedemann \cite{AkramovWiedemann} for the barotropic compressible Euler equations. In his paper \cite{Feireisl}, Feireisl {shows} that weak solutions constructed by convex integration for this problem  cannot satisfy the entropy rate criterion. {We} show these weak solutions cannot satisfy LAAP as well. Section 4 considers the two-dimensional Riemann problem for {an isentropic} fluid as formulated by Chiodaroli and Kreml \cite{ChioKr}. In \cite{ChioKr} the authors produced initial data for which weak solutions constructed by convex integration are preferred according to the entropy rate criterion to the classical two shock solution. {We apply LAAP to give conditions which ensure that for the same data the two shock solution is preferred to the convex integration solutions.}  In Section 5 for the {special} case $p(\rho) = \rho^2$ we extend the local results of Section 4 globally in that, {without further conditions}, the two shock solution is preferred to the convex integration solutions, whenever the convex integration solutions of \cite{ChioKr} exist.


\section{The Least Action Principle}

Here we present a brief explanation of the underlying mechanical
motivation for the least action principle.

Consider the motion of a particle of constant mass $m$ whose position in $\R^3$
at time $t$ is $\mq(t)$.  The particle is acted upon by a potential
$U(|\mq|)$. Form the Lagrangian
\[
  L = \frac12 m |\dot{\mq}(t)|^2 - U(|\mq(t)|),
\]  
so that the action is given by
\[
  A(\mq) = \int_{t_0}^{t_1} \left[\frac12 m |\dot{\mq}(t)|^2 - U(|\mq(t)|)\right]
  \dt.
\]  
Set $|\dot{\mq}(t)|= {\displaystyle \frac{ds}{dt}}$, where $s$ is the
  arc-length along a trajectory. The first expression in the action
  can be written as
  \[
    \begin{eqalign}
\int_{t_0}^{t_1}  \frac12 m |\dot{\mq}(t)|^2 \dt = \frac{m}{2} &
\int_{t_0}^{t_1}\frac{ds}{dt} \frac{ds}{dt} \dt\\
=~ \frac{m}{2} & \int_0^{d_{tot}} \frac{ds}{dt} \, \hbox{d} s,
\end{eqalign}
\]
where $d_{tot}$ is the total distance travelled by the particle. Then
we have
\[
 \begin{eqalign} 
   \int_{t_0}^{t_1}  \frac12 m & |\dot{\mq}(t)|^2 \dt =
   \frac{m}{2} \frac{d_{tot}}{d_{tot}} \int_0^{d_{tot}} \frac{ds}{dt}
   \,
   \hbox{d} s\\
   &~= \frac{m}{2} d_{tot}\ (\hbox{average speed}).
 \end{eqalign}
\]


The expression $\frac{m}{2} d_{tot} (\hbox{average speed})$ is the
cost or effort of moving the particle in the absence of an applied
force. Notice that this is exactly the definition of action given by
Maupertuis in his classic treatise \cite{Maup}.

\begin{figure}[htb]
\centerline{
\resizebox{7cm}{!}{
\begin{tikzpicture}[ 
  line cap=round,
  line join=round,
  >=Triangle,
  myaxis/.style={->,thick}]

\draw[myaxis] (0,0) -- (5.00,0) node[below left = 0.6mm] {$|\mb{q}|$};      
\draw[myaxis] (0,0) -- (0,5) node[below left= 1mm] {$U$};

\draw[smooth, domain=1:4] plot[smooth] (\x, {1.2+0.7*(\x-2.5)*(\x-2.5)});

\end{tikzpicture}
}}
\caption{\label{U} A quadratic potential $U$. }
\end{figure}
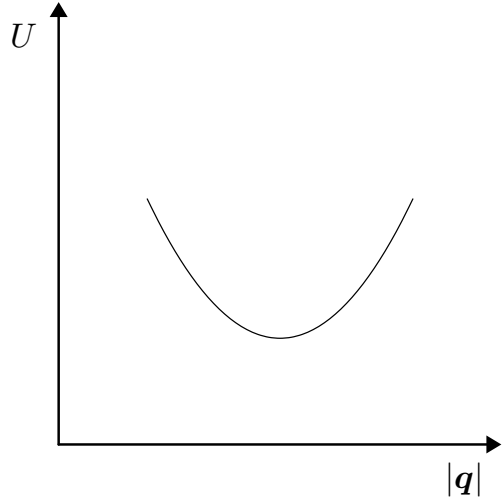

The second term in this action is
\[
  - \int_{t_0}^{t_1} U(|\mq (t)|) \dt.
\]
Imagine our particle moves in a parabolic potential well, as in Figure \ref{U}.

It will require more effort from a deeper potential well than a
shallower one. Thus we wish to penalize for more effort and hence the
minus sign.

We see that taken together the two integrals represent the total cost
or effort of moving the particle in the time interval
$[t_0,t_1]$. When generalized to the barotropic fluid, the action
represents the total cost or effort  to move the fluid in the time
interval $[t_0,t_1]$. Hence the classical least action principle says
that material motion should minimize the total cost or effort of the
motion.


From the above definitions of the action one can formally recover the
equations of particle motion and the inviscid Euler equations of both
the barotropic compressible fluid and an incompressible fluid by
finding necessary conditions for a stationary point. There is a small
proviso here. The computations used to find the Euler equations for
fluid flow are best done in Lagrangian (material) coordinates. This
computation can be done in Eulerian coordinates, but the Lagrangian
must be adjusted with the use of Lagrange multipliers and has been the
subject of considerable discussion in the literature (see Bateman
\cite{Bateman} and Seliger and Whitham \cite{SeligerWhitham}).


\section{The Dafermos oscillator}

In \cite{Dafermos2}, C.~M.~Dafermos considered the second order
nonlinear oscillator
\begin{equation}\label{daf}
  \ddot{x}=g(x,\dot{x}),
\end{equation}
where
\begin{equation}\label{eqg}
  g(x,\dot{x}) = \begin{cases}
    \frac{\dot{x}^2-x^2}{2x},  &$(x-1)^2+ \dot{x}^2 \geq 1, \; x> 0$,\\
    1-x,                      & $(x-1)^2+ \dot{x}^2 <1$.\\
  \end{cases}
\end{equation}

For initial data $(x(0),\dot{x}(0))$ in the strict right half plane
exterior to the unit circle centered at $x=1$, $\dot{x}=0$, trajectories
move clockwise on circles of radius
${\displaystyle c_0 = \frac{x(0)^2+\dot{x}(0)^2}{2x_0} > 1}$ until
they reach the origin $x=\dot{x}=0$. At the origin, $g$ is not
defined, and the trajectories may switch to another circle of radius
$c>1$, thus imposing a jump condition in the acceleration $\ddot{x}$
across $x=0$, $\dot{x}=0$. Dafermos applied both the method of
limiting artificial viscosity and the entropy rate criterion to assert
that at the origin according to these criteria solutions switch to the
unit circle centered at $x=1$, $\dot{x}=0$. This issue was pursued
again in our earlier paper \cite{GGKS}.


Here we re-examine the Dafermos oscillator problem in the light of 
LAAP. The Lagrangian for (\ref{daf}) is
\begin{equation}\label{LagDaf}
  L =  \frac{\dot{x}^2}{2x} - \frac{x}{2},  \; \; (x-1)^2+ \dot{x}^2 \geq
  1, \; x> 0.
\end{equation}

For a trajectory exiting but not touching the origin on a circle with
$c>1$ centered at $x=c$, $\dot{x}=0$, we have that
\[
  L=c-x,
\]
since
\[
  x=c \cos(-t+\theta_0) +c, \; \; \dot{x} = c \sin (-t + \theta_0).
\]
Choose $\theta_0=\pi$ so that our trajectory is exiting the
origin. Then
\[
  L=-c \cos(-t+\pi),
\]
and the action on the time interval $[0,t_1]$ is
\[
 \begin{eqalign} 
  A =~& -c \int_0^{t_1} \cos(-t +\pi) \dt\\   
  =~& -c (-\sin(-t+\pi)) |_0^{t_1} \\
  =~& c \sin(-t_1 + \pi).
\end{eqalign}
\]
For $0 < t_1 < \pi$, the action is minimized by choosing the smallest
available value of $c$. But the best choice is $c=1$. Thus we see that  
LAAP produces the same admissible trajectory as both the limiting
artificial viscosity and the entropy rate criteria.


\section{The barotropic compressible Euler system}

In this section we consider the barotropic compressible Euler system
\begin{eqalignno}  
  &~\partial_t \rho + \hbox{div}_x\, \mb{m} = 0, \label{Euler-mass}\\
  &~\partial_t \mb{m} + \div_x{\left( \frac{\mb{m} \otimes \mb{m}}{\rho}
  \right)} + \nabla_x p(\rho)=0,\label{Euler-mom}
\end{eqalignno}
with unknown density $\rho$, momentum $\mb{m}$ and constitutively
determined pressure $p(\rho)$. The fluid velocity is
$\mb{v}={\displaystyle \mb{m}/\rho}$ and all dependent variables
depend on position $\mb{x}$ and time $t$,
$(t,\mb{x}) \in [0, \infty) \times \R^2$. We are interested in the
initial value problem, where $(\rho, \mb{m})$ satisfy the initial
condition
\begin{equation}\label{Euler-init}
  (\rho, \mb{m})(0,\cdot)= (\rho_0, \mb{m}_0).
\end{equation}

A pair $(\rho, \mb{m})$ is a weak solution of the initial value
problem for the barotropic compressible Euler equations if
(\ref{Euler-mass})--(\ref{Euler-init}) are satisfied in the sense of
distributions.


A weak solution  $(\rho, \mb{m})$ of
(\ref{Euler-mass})--(\ref{Euler-init}) in an open set $\Omega\subset \mathbb{R}^n$ is said to satisfy the
\emph{energy--entropy inequality} if
\begin{equation}\label{EEineq}
\partial_t \left( \frac{|\mb{m}|^2}{2\rho} + \rho \epsilon(\rho)
\right) + \div_x \left[ \left( \frac{|\mb{m}|^2}{2\rho} + \rho \epsilon(\rho)
    + p(\rho) \right) \frac{\mb{m}}{\rho} \right] \leq 0
\end{equation}
in the sense of distributions, where $\epsilon(\rho)$ denotes the specific internal energy per unit mass, {i.e., $\epsilon'(\rho) = \frac{p(\rho)}{\rho^2}$}. 

A weak solution  $(\rho, \mb{m})$ of
(\ref{Euler-mass})--(\ref{Euler-init}) is said to satisfy the 
\emph{entropy rate admissibility criterion} if there is no other weak solution
$(\overline{\rho}, \overline{\mb{m}})$ with the property that for some
$\tau \in [0,\infty)$, $(\rho, \mb{m})= (\overline{\rho},
\overline{\mb{m}})$ on $[0,\tau] \times \Omega$ and
\[
  D[\overline{\rho}, \overline{\mb{m}}](\tau) < D [\rho,\mb{m}](\tau),
\]
where
\begin{equation}\label{ER}
  D [\rho,\mb{m}](t) \stackrel{\cdot}{=} \frac{d^+}{dt} \int_\Omega  \left(
    \frac{|\mb{m}|^2}{2\rho} + \rho \epsilon(\rho) \right) \dx
\end{equation}
{with a corresponding expression for $D[\overline{\rho}, \overline{\mb{m}}](t)$.}

We conclude that the entropy rate criterion prefers the weak solution for
which the total energy decreases with maximal rate.

In this regard, E.~Feireisl \cite{Feireisl} has shown that weak
solutions satisfying the energy--entropy inequality (\ref{EEineq})
which are generated by convex integration do not comply with the
entropy rate criterion. In particular, he shows that for any such
solution there is another such solution which has a larger
dissipation rate.

In this section we use a result by Akramov and Wiedemann
\cite{AkramovWiedemann} to illustrate Feireisl's argument and, more
importantly, how LAAP gives a result similar to Feireisl's. {The solutions of \cite{AkramovWiedemann}  are also constructed by convex integration.}

\begin{theo}[\cite{AkramovWiedemann}] \label{AW1}
Let $n \geq 2$, $\Omega \subset \R^n$ a bounded open set, $T>0$ and $
\Omega' \supset \supset \Omega$ locally Lipschitz. {For a positive constant $\overline{\rho}$,} assume that
$\rho_0 \in C^1(\R^n)$ is a positive function satisfying $\rho_0
(\mb{x}) = \overline{\rho} > 0$ for $\mb{x} \in \R^n \backslash
\Omega$ and the pressure $p \in C^1(\R^n)$ satisfies
\[
  \int_\Omega p(\rho_0(\mb{x})) \, \hbox{{\em d}}\mb{x} =
    p(\overline{\rho}) \,\hbox{{\em meas}}\, \Omega.
\]

Then there exists a bounded initial momentum $\mb{m}_0$ with
$\hbox{{\em supp}}\, (\mb{m}_0) \subset \Omega'$ for which there are
infinitely many solutions $(\rho, \mb{m})$ of
(\ref{Euler-mass})--(\ref{Euler-init}) with density
$\rho(\mb{x})=\rho_0(\mb{x})$.

Moreover, for the obtained weak solution $\mb{m}$ satisfies
\[
  |\mb{m}(t,\mb{x})|^2 = \rho_0(\mb{x})\chi(t) \1_{\Omega'} \hbox{
    a.e. in } [0,T) \times \R^n,
\]
\[
  |\mb{m}_0(\mb{x})|^2 = \rho_0(\mb{x})\chi(0) \1_{\Omega'} \hbox{
    a.e. in } \R^n,
\]
for some smooth function $\chi: \R \mapsto \R$. 
\end{theo}

\begin{theo}[\cite{AkramovWiedemann}] \label{AW2}
  Under the same assumptions as Theorem~\ref{AW1}, there
  exists a maximal time $\overline{T}>0$ such that the weak solutions
  $(\rho, \mb{m})$ satisfy the energy--entropy inequality
  {\em (\ref{EEineq})} in the sense of distributions.
\end{theo}

\begin{cor}[\cite{AkramovWiedemann}]\label{AW3}
  Let $n \geq 2$ and $\Omega \subset \R^n$ be a non-empty bounded open
  set. Assume that $\rho_0 \in C^1(\R^n)$ satisfies
  $\rho_0(\mb{x}) >0 $ for every $\mb{x} \in \R^n$,
  $\rho_0 (\mb{x}) = \overline{\rho} > 0$ for
  $\mb{x} \in \R^n \backslash \Omega$. Let   
  $p \in C^1(\R^n)$ be a given function satisfying
\[
  \int_\Omega p(\rho_0(\mb{x})) \, \hbox{{\em d}}\mb{x} =
    p(\overline{\rho}) \,\hbox{{\em meas}}\, \Omega.
\]
Then there exist $\Omega' \supset \Omega$, $\mb{m}_0$ and a positive
time $\overline{T}$ such that
$\hbox{{\em supp}}\, \mb{m}_0 \subset \Omega'$,
$\hbox{{\em div}}\, \mb{m}_0 =0$, for which there exist infinitely
many $\mb{m}$ such that
$\hbox{{\em supp}}\, \mb{m} (t, \cdot) \subset \Omega'$ for
$t \in [0, \overline{T})$ and $(\rho, \mb{m})$ is a weak solution of
{\em (\ref{Euler-mass})--(\ref{Euler-init})} which satisfies {\em
  (\ref{EEineq})} in the sense of distributions on $[0,\overline{T})
\times \R^n$  with $ \rho(t,\mb{x})  = \rho_0(\mb{x})
\1_{[0,\overline{T})}(t) \in C^1([0,\overline{T}) \times \R^n)$.
\end{cor}

{Substitution of the results given in Theorems~\ref{AW1}, \ref{AW2}
gives
\[
  \int_{\Omega'} \left[\frac{|\mb{m}|^2}{2 \rho} + \rho \epsilon
    (\rho) \right] \dx = \frac12 \chi(t) \, \hbox{meas} \,\Omega' +
  \int_{\Omega'} \rho_0 \epsilon (\rho_0) \dx \hbox{ on }
  [0,\overline{T}),
\]
and hence 
\begin{equation}\label{AWen}
  \frac{d^+}{dt} \int_{\Omega'} \left[\frac{|\mb{m}|^2}{2 \rho} + \rho \epsilon
    (\rho) \right] \dx = \frac12 \chi'(t) \, \hbox{meas} \,\Omega'.
\end{equation}}
Akramov and Wiedemann provide a sufficient condition for the
non-unique weak solutions to satisfy the energy--entropy inequality; namely, $\chi(t)$ should be chosen to satisfy
\begin{equation}\label{chi1}
  \chi'(t) \leq -C_1 \chi^{1/2}(t) - C_2 \chi^{3/2}(t),
\end{equation}
and
\begin{equation}\label{chi2}
  \chi(t) >  n \lambda(t) \hbox{ on } [0, T],
\end{equation}
where ${\lambda \geq 0}$ is determined by their choice of sub-solutions
and the {positive} constants $C_1$, $C_2$ are determined by estimates on $\nabla
\rho_0(\mb{x})$.  Clearly $\chi (t)$ is dominated by solutions of
\[
  \chi_1'(t) =-C_1 \chi_1^{1/2}(t), \; \; \chi_1(0)= \chi(0), 
\]
which are given by
\[
  \chi_1(t) = \left(\chi^{1/2}(0) - \frac12 C_1 t\right)^2,
\]
and we see that $\chi_1(2\chi^{1/2}(0)/C_1)=0$. Hence $\chi(t)$ satisfying
(\ref{chi1}) is monotone decreasing and vanishes for some
$\overline{T} < {\displaystyle \frac{2\chi^{1/2}(0)}{C_1}}$. 

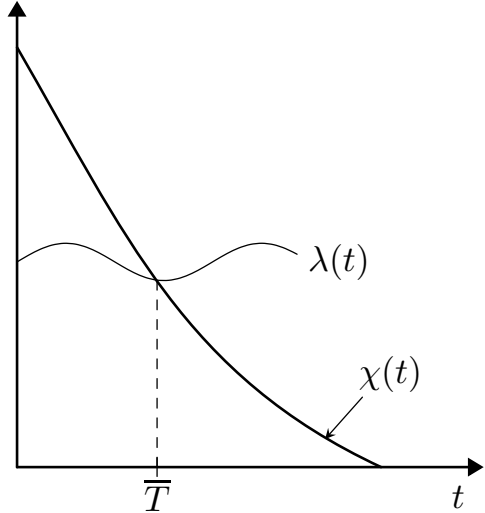
\begin{figure}[htb]
\centerline{
\resizebox{7cm}{!}{
\begin{tikzpicture}[ 
  line cap=round,
  line join=round,
  >=Triangle,
  myaxis/.style={->,thick}]

\draw[myaxis] (0,0) -- (5.00,0) node[below left = 0.6mm] {$t$};      
\draw[myaxis] (0,0) -- (0,5) node[below left= 1mm] {};

\draw (1.5, -0.1) -- (1.5,0.1);
\node at (1.5,-0.3) {$\overline{T}$};

\draw[thick] (0.0, 4.5) to [out=-60, in=155] (3.9,0.0);
\draw[dashed] (1.5,0.0) -- (1.5,2.05);

\draw [-stealth] (3.7, 0.75) -- (3.3, 0.3);
\node at (4.00,1.00)  {$\chi (t)$};

\node at (3.45,2.20)  {$\lambda (t)$};

\draw[smooth, domain=0:3] plot[smooth] (\x, {2.2+0.2*sin(3*(\x r))});

\end{tikzpicture}
}}
\caption{\label{chi} Behavior of $\chi (t)$. }
\end{figure}

Inspection of a typical graph of $\chi(t)$, as in Figure \ref{chi}, shows that for any choice
of $\chi(t)$ satisfying (\ref{chi1}--\ref{chi2}), we could always
choose larger constants $C_1, C_2$ and a possible shorter time
interval so that the new choice of $\chi(t)$ will decrease more
rapidly.  Thus consistent with the earlier quoted result of Feireisl,
there  is no solution given by Akramov and Wiedemann that satisfies the
entropy rate admissibility criterion.

Now let us compute the action for the Akramov and Wiedemann
solution, obviously given by
\[
  \begin{eqalign}
    \int_0^{\overline{T}} & \int_{\Omega'} \left[\frac{|\mb{m}|^2}{2
        \rho} - \rho \epsilon (\rho) \right] \dx \dt\\ &~=
    \frac12 \left( \int_0^{\overline{T}} \chi(t) \dt \right) \,
    \hbox{meas} \, \Omega' - \overline{T} \int_{\Omega'} \rho_0
    \epsilon (\rho_0) \dx. 
 \end{eqalign}   
\]
Hence the action will be minimized by the smallest choice of $L^1([0,
\overline{T}))$ norm of $\chi$, i.e.~$\chi=0$ a.e.~{on $[0,T]$}. But by
Theorem~\ref{AW1}, this gives $\mb{m} = 0$ and hence this choice
will not yield a solution to the initial value problem unless
$\mb{m}_0 \equiv 0$ a.e. Thus LAAP yields a result consistent
with the previously stated conclusion from the entropy rate admissibility criterion,
i.e.~there is no solution given by Akramov and Wiedemann that satisfies
LAAP.   

\section{Convex integration and the Riemann problem}

In this section we consider admissibility of weak solutions to the
Riemann initial value problem and in particular compare the entropy rate admissibility criterion with that of LAAP. 
{We suppose
familiarity with the presentation by Chiodaroli and Kreml \cite{ChioKr}  and the more recent one
by Markfelder \cite{Markfelder}.}

{We consider the isentropic form of the compressible Euler
system (\ref{Euler-mass})--(\ref{Euler-init}) in which $p(\rho) = \rho^\gamma$, $\gamma \geq 1$, and prescribe Riemann initial data to
be}
\begin{equation}\label{Rie}
  (\rho_{{0}}, \mb{m}_{{0}}) = \begin{cases}
    (\rho_-, \mb{m}_-), & $x_2<0$, $-\infty < x_1 < \infty$, \\
    (\rho_+, \mb{m}_+), & $x_2>0$, $-\infty < x_1 < \infty$,
  \end{cases}
\end{equation}  
{and $\rho_\pm, \mb{m}_\pm$ are constants.} {We now summarise definitions and concepts
required for the discussion of weak solutions constructed in \cite{ChioKr} using convex integration.
We use the general relation $\textbf{m}=\rho \textbf{v}$ to  replace the momentum $\textbf{m}$.}

{Let $S^{2\times 2}_{0}$ denote the set of $2\time 2$ symmetric matrices with zero trace.}

{(i) A \textit{fan sub-solution} is a triple $(\bar{\rho},\,\bar{\textbf{v}},\,\bar{\textbf{u}}): \mathbb{R}^{2}\times (0,\,\infty)\rightarrow \mathbb{R}^{+}\times\mathbb{R}^{2}\times S^{2\times 2}_{0}$ of piecewise constant functions such that
\begin{equation}\label{fan}
  (\overline{\rho},\overline{\mb{v}}, \overline{\mb{u}}) = (\rho_-,\mb{v}_-, \mb{u}_-) \1_{P_-} + (\rho_1,\mb{v}_1,
  \mb{u}_1) \1_{P_1} +  (\rho_+,\mb{v}_+, \mb{u}_+) \1_{P_+},
\end{equation}
where $(\rho_{1},\textbf{v}_{1},\textbf{u}_{1})$ are constants with $\rho_{1}>0$ which   for positive constant $C$ satisfy
\begin{equation}
\label{condii} 
\textbf{v}_{1}\otimes \textbf{v}_{1}-\textbf{u}_{1}<\frac{C}{2}Id.
\end{equation}
We adopt the general notation,
\begin{equation}\label{compu}
{u}_{ij}= {v}_i {v}_j - \frac12 |{ \mb{v}}|^2 \, \hbox{Id},
\end{equation}
and let}
\begin{equation}\label{Ps}
  \begin{eqalign}
    &~P_- = \{ (x_1, x_2,t), \; x_2 < \nu_- t,\; -\infty < x_1 < \infty, \; t>0
    \},\\
    &~P_1 = \{ (x_1, x_2,t), \; \nu_-t < x_2 < \nu_+ t,\;  -\infty < x_1
    < \infty, \; t>0 \},\\
    &~P_+ = \{ (x_1, x_2,t), \; \nu_+ t < x_2 ,\;  -\infty < x_1 <
    \infty, \; t>0\},
  \end{eqalign}
\end{equation}  
{for some numbers $\nu_-,\nu_+$, $0< \nu_-<\nu_+$. The shock speeds $\nu_-,\nu_+$ depend upon $\overline{\rho}$, and we write $\nu_-(\overline{\rho}),\nu_+(\overline{\rho})$ to emphasize this dependence.}\\

{(ii) With the constant $C$ from \eqref{condii} the triple $(\bar{\rho},\,\bar{\textbf{v}},\,\bar{\textbf{u}})$ solves}
\begin{equation}\label{condiii}
  \partial_t \overline{\rho} + \hbox{div}_{\mb{x}} \, (\overline{\rho}\
  \overline{\mb{u}} ) = 0,
\end{equation}
\begin{equation}\label{condiiia}
    \partial_t (\overline{\rho}\ \overline{\mb{v}}) + \hbox{div}_{\mb{x}} (\overline{\rho}\
    \overline{\mb{u}}) + \nabla_{\mb{x}} \Big[ p(\overline{\rho}) +
      \frac12 C\rho_1 \1_{P_1} 
    +  \overline{\rho} |\overline{\mb{v}}|^2 \1_{P_- \cup P_+} \Big] =
    0
\end{equation}  
in the sense of distributions. 

Notice that in $P_-, P_+$, $(\overline{\rho}, \overline{\mb{v}},
\overline{\mb{u}})$ is just the usual weak solution of (\ref{Euler-mass})--(\ref{Euler-init}) but in $P_1$ we have
\begin{equation}\label{eqp1}
    \partial_t (\overline{\rho} \overline{\mb{v}}) + \hbox{div}_{\mb{x}} (\overline{\rho}
    \overline{\mb{u}}) + \nabla_{\mb{x}} \left( p(\overline{\rho}) +
      \frac12 C\rho_1\right) =0. 
  \end{equation}

The crucial lemma that provides existence of non-unique weak solutions
to (\ref{Euler-mass}), (\ref{Euler-mom}), (\ref{Rie}) is Lemma 3.2
of \cite{ChioDeLKr} which originates in the work of Tartar \cite{Tartar}.


\begin{lemma}[Lemma 3.2 of \cite{ChioDeLKr}] \label{CK32}
  Let $(\tilde{\mb{v}}, \tilde{\mb{u}}) \in \R^2 \times S_0^{2\times 2}$
  and $C>0$ be such that $\tilde{\mb{v}} \otimes \tilde{\mb{v}} -
      \tilde{\mb{u}} < {\displaystyle \frac{C}{2} \hbox{{\em Id}}}$, {as in \eqref{condii}}.
        For any $\Omega \subset \R^+ \times \R^2$ there are infinitely
        many maps $(\mb{v}, \mb{u})$ with the following properties:
\begin{enumerate}[(i)]       
\item $\mb{v}$ and $\mb{u}$ vanish identically outside $\Omega$;
\item $\hbox{{\em div}}_{\mb{x}} \mb{v} =0$ and $\partial_t \mb{v} +
  \hbox{{\em div}}_{\mb{x}} \mb{u} =0$ in the sense of distributions;
\item $(\tilde{\mb{v}}+ \mb{v}) \otimes (\tilde{\mb{v}}+ \mb{v})-
  (\tilde{\mb{u}}+ \mb{u}) = {\displaystyle \frac{C}{2}} \,
  \hbox{{\em Id}}$ a.e. in $\Omega$.
\end{enumerate}  
\end{lemma}

The key observation in applying Lemma~\ref{CK32} is to add
$(\mb{v}, \mb{u})$ of Lemma~\ref{CK32} to the fan sub-solution, i.e.~take
{$\tilde{\mb{v}}= \mb{v}_1$, $\tilde{\mb{u}}=\mb{u}_{1}$,} $\Omega=
P_1$. Obviously,
\[
  \partial_t (\rho_1(\mb{v}+\mb{v_1}))=0
\]
and also
\begin{equation}\label{eq49}
    \partial_t (\rho_1(\mb{v} + \mb{v}_1) + \hbox{div}_{\mb{x}} (\rho_1
    (\mb{u} + \mb{u}_1)) + \nabla_{\mb{x}} \left[ p(\rho_1)
      + \rho_1\left( ( \mb{v}+\mb{v}_1) \otimes ( \mb{v}+\mb{v}_1)
        - (\mb{u}+\mb{u}_1) \right)\right]=0, 
  \end{equation}
where we used (iii) of Lemma~\ref{CK32}.

Next observe exact cancellation of the second and last terms in
(\ref{eq49}). In consequence, $(\rho_1, \mb{v}+\mb{v}_1)$
satisfy (\ref{Euler-mass}), (\ref{Euler-mom}) in $P_1$ and we have
produced an infinite number of weak solutions to the Riemann initial
value problem. Further, (iii) of Lemma~\ref{CK32} implies
\begin{equation}\label{eq410}
  |\mb{v}+ \mb{v}_1|^2\1_{P_1} = C.
\end{equation}
 {Now the shock speeds are given by $\overline{\rho} = \rho_1$, i.e., $\nu_-(\rho_1),\nu_+(\rho_1)$ (see \cite[(4.46), (4.47), (4.51), (4.52)]{ChioKr} for the relevant formulas).}


The constant $C$ given in \eqref{eq410} obeys the following identity \cite[p.~1043]{ChioKr}:
\begin{equation}\label{eq411new}
C=\beta^2(\rho_1)+\epsilon_1(\rho_1)+\epsilon_2,
\end{equation}
where $\beta(\rho_1)$ is given by \cite[(4.48), (4.53)]{ChioKr}, $\epsilon_1(\rho_1)$ is given by \cite[(4.49), (4.50), (4.54)]{ChioKr} and $\epsilon_2$ is any positive constant satisfying \cite[(4.82), (4.83)]{ChioKr}. In particular, $\beta$, $\epsilon_1$ are smooth functions of $\rho_1$, $\beta^2(\rho_m) = v_m^2$, $\epsilon_1(\rho_m)=0$ and $\epsilon_2>0$ on some interval {\begin{equation}\label{Ilabel}I=[\rho^\ast,\rho_m].\end{equation}}

The values $\rho_m$, $v_m$ refer to density, velocity of the two shock solution to the Riemann initial value problem
\begin{equation}\label{eq411}
  \begin{cases}
    \rho= \rho_-, \; v_a=0, \; v_b = v_{b-}, & $x_2<\nu_-(\rho_m) t$,\\
    \rho= \rho_+, \; v_a=0, \; v_b = v_{b+}, & $x_2>\nu_+(\rho_m) t$,\\
    \rho= \rho_m, \; v_a=0, \; v_b = v_m, & $\nu_-(\rho_m) t<x_2<\nu_+(\rho_m) t$,
  \end{cases}
\end{equation}  
{where $\mb{v} =(v_a, v_b)$, and} $\nu_-(\rho_m), \nu_+(\rho_m)$,  $\rho_-, \rho_+, \rho_m$, $v_{b-}, v_{b+}, v_m$ are consistent with the Rankine-Hugoniot conditions.

%
%
%
%
%
%
%
%


 We now define the action of both the
Chiodaroli--Kreml weak solutions and the two shock solution to be:
\[
  \begin{eqalign}
  \int_{-L_3}^{L_3} \int_0^T \int_{\ell_1}^{\ell_2}\mathcal{L}\
  \hbox{d}x_2 \dt \hbox{d}x_1 ,
\end{eqalign}  
\]
where $\mathcal{L}$ denotes the {respective} Lagrangian, $\ell_1<\min\{\nu_-(\rho_1) : \rho_1 \in I\}T$ and $\ell_2>\max\{\nu_+(\rho_1) : \rho_1 \in I\}T$ and $T>0$ is fixed. Set

\begin{equation}\label{eq412new}
D(\rho_1) = \int_{-L_3}^{L_3} \int_0^T \int_{\ell_1}^{\ell_2} L_{diff}(\rho_1) \ \mathrm{d}x_2 \dt \ \mathrm{d}x_1 ,
\end{equation}
where \[
  L_{diff} = \hbox{Lagrangian of two shock solution } 
- \hbox{Lagrangian of convex integration solution}.
\]
$D(\rho_1)$ is a continuous function of $\rho_1$. {Since $\rho_{1}\in I$ we may set $\rho_{1}=\rho_{m}$ to obtain}
$$D(\rho_m) = \int_{-L_3}^{L_3} \int_0^T \int_{\nu_-(\rho_m) t}^{\nu_+(\rho_m) t} L_{diff}(\rho_m) \ \mathrm{d}x_2 \dt \ \mathrm{d}x_1.$$

We have $$L_{diff}(\rho_1) = \frac12 \rho_m v_m^2 - \rho_m \epsilon(\rho_m) - \frac12 \rho_1
  C + \rho_1 \epsilon (\rho_1)$$
and by \eqref{eq411new}, {when $\rho_1=\rho_m$,} \begin{equation}\label{Ldiff0}L_{diff}(\rho_m)=-\frac{1}{2} \rho_m \epsilon_2<0.\end{equation}
We thus see in a {sufficiently} small left neighborhood of $\rho_m$ that $D(\rho_1)<0$, and we can state our main result.

\begin{theo}\label{ThMSnew}
{Let $\epsilon_{2}>0$ and a left interval $I$ of $\rho_{m}$ be fixed so that} for $\rho_1$ contained in this interval a fan subsolution is defined. Via Theorem 1 of \cite{ChioKr} there is an infinite number of convex integration weak solutions to the Riemann problem. 
{Since the action on the interval $I$} of the two shock solution is lower than any of these convex integration solutions, LAAP rules out the admissibility of {the} convex integration solutions and allows for the possible admissibility of the two shock solution. Furthermore  if there are no other entropic weak solutions other than the convex integration solutions and the two shock solution, the two shock solution is the unique LAAP solution.
\end{theo}

We note that this result compares {on the interval $I$} the two shock solution and convex integration solutions where they both satisfy the same initial conditions, behavior as $x_2 \to \pm \infty$ and same periodicity in $x_1$. {An identical} argument will apply to the periodic solutions in $x_2$ constructed in Section 6 of \cite{ChioKr}, where now the spatial domain is a unit square and $T$ is chosen sufficiently small.

Theorem  \ref{ThMSnew} is obviously the exact opposite of the results in \cite{ChioKr}, where Riemann data are given such that the entropy rate admissibility criterion rules out the two shock solution.

\section{Global result for {$p(\rho) = \rho^2$}}

The following main result gives a stronger, global version of Theorem \ref{ThMSnew}. We consider the pressure law $p(\rho) = \rho^\gamma$ with $\gamma=2$, when the specific internal energy is $\rho$. Numerical evidence supports the conjecture that an analogous result holds for all $\gamma \in
[1,3]$.

Recall the necessary condition for the existence of the two shock solution, Equation (2.47) in \cite{ChioKr}:
\begin{equation}\label{247}
  S(v_{b,-},v_{b,+},\rho_+,\rho_-) := (v_{b,-}-v_{b,+})^2\rho_+\rho_-
  -(\rho_+ - \rho_-)(p(\rho_+)-p(\rho_-)) > 0.
\end{equation}  

\begin{theo} \label{global} Let $\gamma=2$.  For any Riemann data
  satisfying \eqref{247} and any $\epsilon_2>0$ in \eqref{eq411new},  the action of the two shock solution is lower than the action of any of the convex integration solutions constructed in Theorem 1 of \cite{ChioKr}. LAAP therefore rules out the admissibility of these convex integration solutions and allows for the possible admissibility of the two shock solution. Furthermore  if there are no other entropic weak solutions other than the convex integration solutions and the two shock solution, the two shock solution is the unique LAAP solution.
\end{theo}
The proof of the theorem involves detailed computations, which are performed using the Maple worksheet \cite{MapleWS}. As in \cite{ChioKr}, we restrict
ourselves to the case of $R = \rho_- -\rho_+>0$; the analysis of the other cases is similar. 

Recall from above that 
\[
  L_{diff} = \hbox{Lagrangian of two shock solution } 
- \hbox{Lagrangian of convex integration solution}.
\]
We express $L_{diff}$ in terms of the variable
  \[
    w= \sqrt{\frac{\rho_1-\rho_-}{\rho_1-\rho_+}}.  
\]
The result is that
  \[
    L_{diff}(\rho_1) = \frac{P(w,\epsilon_2)}{Q(w)}, 
  \]
where the polynomial $P(w,\epsilon_2)$ is of degree $6$ in $w$ and
linear in $\epsilon_2$, while \begin{equation}\label{Qeq}Q(w) = -2 (\rho_- \rho_+)  w \rho_- (1-w) (1+w) (\rho_- - \rho_+ w^2)\end{equation} is a polynomial of degree $5$ in $w$.\footnote{{This is verified in the code \cite{MapleWS}, lines 1--76, the formula for $Q$ in  line 84.}}

The following result shows that it suffices to consider $\epsilon_2=0$.
\begin{prop}\label{e2}
  The coefficient of $\epsilon_2$ in $L_{diff}$ is negative.
\end{prop}
\begin{proof}
  One computes\footnote{{This is verified in the code \cite{MapleWS}, line 76.}}   that this coefficient is equal to
  \[
    -\frac{(\rho_+ w^2 +\rho_-)\sqrt{S(v_{b,-},v_{b,+},\rho_+,\rho_-)}}{
      2(\rho_--\rho_+) w}.
  \]
The result follows, because $w \in (0,1)$ and $R = \rho_--\rho_+>0$.
\end{proof}

In the following  we set $\epsilon_2=0$ and denote $P(w,0)$ by $P(w)$.
We now consider the denominator $Q(w)$.
\begin{prop}\label{q}
  For $w \in (0,1)$, $Q(w)$ is negative.
\end{prop}  
\begin{proof}
The assertion follows from the expression \eqref{Qeq}.\end{proof}
In the next step we simplify the numerator $P(w)$. Two roots of $P(w)$ are given by $\pm b$, where 
\[
 b = \sqrt{\frac{\rho_-}{\rho_+}} > 1.
\]
This allows us to factorize
\[
  P(w)=(b^2-w^2)N(w),
\]
where $N(w)$ is a polynomial of degree $4$, which we now analyse.

$N(w)$ has at least one root in
$[0,1]$, namely
  \[
    w^* = \sqrt{\frac{\rho_m-\rho_-}{\rho_m-\rho_+}},
  \]
where $\rho_m$ is given by (2.48) of \cite{ChioKr}. To see this, note that $w^*$ is a root of
$N(w)$ because, according to Equation \eqref{Ldiff0} with $\epsilon_2=0$, $L_{diff}$ has a zero at $\rho_1=\rho_m$.

We conclude that in order to prove Theorem~\ref{global}, it suffices to show that $N(w)$ is positive for 
$w \in (0,w^*)$. Indeed, recall that $Q(w)$ is negative
for $w \in (0,1)$ and that $b^2-w^2>0$ for $w \in (0,1)$, since
$b>1$. 

Therefore, Theorem~\ref{global} follows, if we show that 
\begin{enumerate}
\item[(1)] $N(0)>0$ and 
\item[(2)] $N(w)$ does not have a root in $(0,w^*)$.
\end{enumerate}
By computation, we observe claim (1):
\begin{prop}\label{coefs}
The coefficients $c_0, \ldots, c_4$ of $N(w) = c_0 + c_1 w + c_2 w^2 + c_3 w^3 + c_4 w^4$ satisfy: $c_0>0$, $c_1=-c_3$, $c_4=-c_0$,
$c_2< 0$. In particular, $N(0)=c_0>0$. 
\end{prop}
\begin{proof}
The properties are true by inspection of the explicit formulae for $c_j$. Specifically,
\begin{align*}    
c_0  &= -c_4 = S(v_{b,-}, v_{b,+},\rho_+, \rho_-)^{3/2} \rho_+
\rho_->0,\\
c_2&=-3\sqrt{S(v_{b,-},v_{b,+},\rho_+,\rho_-)}\rho_-\rho_+
    (\rho_--\rho_+)^3 < 0.
\end{align*}
The expression for $c_1=-c_3$ is longer.
\end{proof}


\begin{proof}[\textit{Proof of Theorem \ref{global}}.] Define the auxiliary 
  polynomial
\[
  M(w,t) = t(c_0 +c_2 w^2 - c_0w^4) +c_1w-c_1w^3, 
\]
which satisfies $M(w,1)= N(w)$.

Proposition~\ref{coefs} implies that $w=0$ and $w=\pm 1 $
cannot be roots of $M(w,t)$ for any $t>0$, because $M(0,t)= tc_0>0 $ and
$M(\pm 1,t)= t c_2<0$.  

If $c_1=0$, by solving the quadratic in $w^2$, it
is clear that for all $t>0$ the polynomial $M(w,t)$ has a single positive root, 
and therefore also in particular for $t=1$. 

Now consider the
(generic) case of $c_1 \neq 0$. The roots of $M(w,0)=0$ are then given by $0, \pm 1$. If $t>0$, a fourth real
root comes from $\pm \infty$, but it  never crosses the lines $w = \pm 1$. The other three roots
admit a regular perturbation expansion for $t \to 0^+$. Using the
notation $\kappa_s(t)$ for the expansion of the root that is equal to
$s$ when $t=0$, we find
\[
\begin{eqalign}  
  \kappa_{-1} & \sim -1 + \frac{c_2}{2c_1} t + O(t^2),\\
  \kappa_{0} & \sim  -\frac{c_1}{c_0} t + O(t^2),\\
  \kappa_{1} & \sim 1 + \frac{c_2}{2c_1} t + O(t^2).
\end{eqalign}
\]
We know that $c_0>0$ and $c_2<0$. If $c_1>0$ this means that for all
$t>0$ the root $\kappa_1$ will be confined to the interval $(0,1)$, while the
other two roots will be outside this interval for as long as they exist
(as the straight lines $w=0,\pm 1$ cannot be crossed), so we have a
unique root in $[0,1]$. Similarly, if $c_1<0$, for all $t>0$ the root
$\kappa_0$ will be in the interval $(0,1)$, while the other two roots
will be outside this interval, again giving us a unique root in
$[0,1]$.

Claim (2) and, hence, Theorem  \ref{global} follow.
\end{proof}


\section*{Acknowledgement}
The authors are grateful for support from the ICMS Research in Groups scheme. We thank C.~M.~Dafermos, O.~Kreml and S.~Markfelder for their valuable remarks.

\section*{Data availability statement}
The authors declare that the data supporting the findings of this study are available within the paper.


\end{document}